\documentclass{amsart}

\newtheorem{theorem}{Theorem}
\newcommand{\bt}{\begin{theorem}}
\newcommand{\et}{\end{theorem}}
\newtheorem{lemma}{Lemma}
\newcommand{\bl}{\begin{lemma}}
\newcommand{\el}{\end{lemma}}
\newtheorem{corollary}{Corollary}
\newcommand{\bc}{\begin{corollary}}
\newcommand{\ec}{\end{corollary}}
\newtheorem{problem}{Problem}
\newcommand{\bprob}{\begin{problem}}
\newcommand{\eprob}{\end{problem}}

\newtheorem{example}{Example}
\newcommand{\bex}{\begin{example}}
\newcommand{\eex}{\end{example}}

\newcommand{\beq}{\begin{equation}}
\newcommand{\eeq}{\end{equation}}
\newcommand{\benum}{\begin{enumerate}}
\newcommand{\eenum}{\end{enumerate}}

\newcommand{\Q}{\ensuremath{\mathbf Q}}

\newcommand{\C}{\ensuremath{\mathbf C}}

\newcommand{\bsmallmat}{\left(\begin{smallmatrix}}
\newcommand{\esmallmat}{\end{smallmatrix}\right)}

\newcommand{\bmat}{\left(\begin{matrix}}
\newcommand{\emat}{\end{matrix}\right)}

\usepackage[all]{xy}
\usepackage{amssymb,latexsym}

\title[Continuity of the roots]{Continuity of the roots of a polynomial}
\author{Melvyn B.  Nathanson}
\address{Department of Mathematics\\Lehman College (CUNY)\\Bronx, NY 10468} 
\email{melvyn.nathanson@lehman.cuny.edu} 

\author{David A. Ross}
\address{Department of Mathematics\\University of Hawaii at Manoa\\Honolulu, HI 96822} 
\email{ross@math.hawaii.edu}

\subjclass[2010]{12D10, 12E05, 12J10, 12L15}

\keywords{Polynomials, roots of polynomials, deformations of polynomial, continuity of roots, elementary methods.}

\date{\today}

\begin{document}

\begin{abstract}
Let $K$ be an algebraically closed field with an absolute value.  
This note gives a simple proof of the classical result 
that the roots of a polynomial with coefficients in $K$ are continuous functions 
of the coefficients of the polynomial. 
\end{abstract}

\maketitle

\section{Continuity of the roots} 
Let $K$ be an algebraically closed field with an absolute value 
(for example, the field \C\ of complex numbers or the algebraic closure 
of the field $\Q_p$ of $p$-adic numbers). 
We denote by $K[z]$ the ring of polynomials with coefficients in $K$.

Let  $(a_i)_{i\in I}$ and $(b_i)_{i\in I}$ be sequences of elements of $K$ and let $\delta > 0$.  
The sequence $(b_i)_{i\in I}$  is a \emph{$\delta$-deformation}\index{deformation} 
of the  sequence $(a_i)_{i\in I}$ if 
\[
|b_i - a_i| < \delta
\]
for all $i \in I$.  

Let $f(z) = \sum_{i=0}^n a_i z^i$ and $g(z) = \sum_{i=0}^n b_i z^i$ be 
polynomials in $K[z]$ of degree $n$. 
 Let $\delta > 0$.  The polynomial  $g(x)$ is a  
\emph{$\delta$-deformation}\index{deformation} of the polynomial $f(z)$ if 
the sequence $(b_0, b_1,\ldots, b_n)$ of coefficients of $g(z)$ is a 
$\delta$-deformation of the  sequence $(a_0, a_1,\ldots, a_n)$ of coefficients of  $f(z)$.

In the algebraically closed field $K$, every polynomial of degree $n$ in $K[z]$ 
has $n$ roots, counting multiplicity.  Let 
$(\zeta_1,\ldots, \zeta_{n})$ be the sequence of roots of $f(z)$ 
and let $(\omega_1,\ldots, \omega_n)$ be the sequence of roots of $g(z)$, 
where a root of multiplicity $m$ appears $m$ times in the sequence.  
A classical theorem on the continuity of roots of a polynomial states that, 
for every $\varepsilon > 0$, 
there exists $\delta = \delta(f,\varepsilon) > 0$ such that if the polynomial $g(x)$ 
is a $\delta$-deformation of the polynomial $f(x)$, then there is a permutation $\sigma \in S_n$  
such that the rearranged sequence of roots 
$\left( \omega_{\sigma(1)},\ldots, \omega_{\sigma(n)} \right)$ 
is an $\varepsilon$-deformation of $(\zeta_1,\ldots, \zeta_{n})$. 
Most proofs  of this result require sophisticated or technically difficult 
mathematics (cf.~\cite{cuck-corb89}--\cite{whit72}).  
There are also simple proofs using nonstandard analysis, for example ~\cite{Lutz, Ross}; see also \cite{Remm}. 
The object of this note is to give a proof 
that uses only high school algebra and the definition of continuity.

The open ball with center $\zeta \in K$  and radius $\varepsilon > 0$  
is the set 
\[
B(\zeta,\varepsilon) = \{z\in K: |z- \zeta| < \varepsilon \}.
\]

\bl                    \label{ContRoot:lemma:ball} 
Let $\zeta_1,\ldots, \zeta_{\ell}$ be distinct elements of the field $K$.   If
\beq                                   \label{ContRoot:epsilon}
0 < \varepsilon \leq \frac{1}{2} \min\left( |\zeta_j - \zeta_k| : 1 \leq j < k \leq \ell \right)  
\eeq
then the  open balls $B(\zeta_1,\varepsilon)$, 
\ldots, $B(\zeta_{\ell},\varepsilon)$ are pairwise disjoint.
\el

\begin{proof}
It suffices to observe that if  $z \in B(\zeta_j,\varepsilon) \cap B(\zeta_k,\varepsilon)$, 
then   
\[
2\varepsilon \leq  |\zeta_j - \zeta_k| \leq|\zeta_j - z| + |z- \zeta_k| < 2\varepsilon. 
\]
which is absurd.       This completes the proof.  
\end{proof}

Let  $\zeta_1,\ldots, \zeta_{\ell}$ be the  distinct roots of the polynomial $f(z)$ 
and let $\mu_j$ be the multiplicity of the root $\zeta_j$ 
for all $j \in \{1,\ldots, \ell \}$. 
Thus, $\sum_{j=1}^{\ell} \mu_j = n$ and 
\[
f(z) = a_n\prod_{j=1}^{\ell} (z - \zeta_j)^{\mu_j}
\]
with $a_n \neq 0$.  
Let $\varepsilon > 0$.  
A polynomial $g(z)$ of degree $n$ is \emph{$\varepsilon$-aligned}\index{aligned} to $f(z)$ 
 if, for all $j \in \{1,\ldots, \ell\}$, the open ball 
$B(\zeta_j,\varepsilon)$ contains at least $\mu_j$ roots of $g(z)$, counted with multiplicity.

If $\varepsilon$ satisfies inequality~\eqref{ContRoot:epsilon}, 
then, by Lemma~\ref{ContRoot:lemma:ball},   the  open balls $B(\zeta_1,\varepsilon)$, 
\ldots, $B(\zeta_{\ell},\varepsilon)$ are pairwise disjoint and 
 a polynomial $g(z)$ of degree $n$ is $\varepsilon$-aligned to $f(z)$ 
 if and only if, for all $j \in \{1,\ldots, \ell\}$, the open ball 
$B(\zeta_j,\varepsilon)$ contains exactly $\mu_j$ roots of $g(z)$, counted by multiplicity.    
Each root of $g(z)$ 
belongs to one and only one of the balls $B(\zeta_j,\varepsilon)$. 

In Section~\ref{ContRoot:section:preliminary} we prove that if $f(z)$ and $g(z)$ 
are polynomials of degree $n$ such that $g(z)$ is  $\varepsilon$-aligned 
to $f(z)$ for all $\varepsilon > 0$, then $g(z) = cf(z)$ for some $c \neq 0$.

The theorem on the continuity of the roots of a polynomial states 
that for every $\varepsilon > 0$ there exists $\delta  > 0$ 
such that if the polynomial $g(z)$ is a $\delta$-deformation 
of $f(z)$, then the polynomial $g(z)$ is $\varepsilon$-aligned to $f(z)$.  
Section~\ref{ContRoot:section:no-nonzero} considers the special 
case of polynomials whose only root is zero 
and Section~\ref{ContRoot:section:nonzero}
considers polynomials that have a nonzero root.

\section{The inverse theorem}                   \label{ContRoot:section:preliminary} 

We begin with a preliminary observation.

\bt
Let $f(z)$ and $g(z)$ be polynomials in $K[z]$ of degree $n$.  
Then $g(z)$ is $\varepsilon$-aligned 
to $f(z)$ for all $\varepsilon > 0$
if and only if $g(z) = cf(z)$ for some $c \neq 0$. 
\et

\begin{proof}
If $g(z) = cf(z)$ for some $c \neq 0$, then $f(z)$ and $g(z)$ have the same roots 
with the same multiplicities and so $g(z)$ is $\varepsilon$-aligned to $f(z)$ 
for all $\varepsilon > 0$. 

Suppose that $g(z)$ is $\varepsilon$-aligned to $f(z)$ 
for all $\varepsilon > 0$. 
Let  $\zeta_1,\ldots, \zeta_{\ell}$ be the  distinct roots of the polynomial $f(z)$ 
and let $\mu_j$ be the multiplicity of $\zeta_j$ for all $j \in \{1,\ldots, \ell\}$.  
Let $\varepsilon$ satisfy inequality~\eqref{ContRoot:epsilon}. 
The open balls $B(\zeta_1, \varepsilon), \ldots, B(\zeta_{\ell},\varepsilon)$ 
are pairwise disjoint.  
For $j \in \{1,\ldots, \ell\}$, let $\omega_1,\ldots, \omega_{\mu_j}$ be the $\mu_j$ 
not necessarily distinct roots of $g(z)$ that are 
in the open ball $B(\zeta_j,\varepsilon)$.  
If $0 < \varepsilon' < \varepsilon$, then $B(\zeta_j,\varepsilon')\subseteq B(\zeta_j,\varepsilon)$.  
Because the ball $B(\zeta_j,\varepsilon')$ also contains exactly $\mu_j$ roots of $g(z)$, 
counted by multiplicity, we must have  
\[
\{ \omega_1,\ldots, \omega_{\mu_j} \} \subseteq B(\zeta_j,\varepsilon')
\]
and so 
\[
|\omega_i-\zeta_j| < \varepsilon'
\]
for all $i \in \{1,\ldots, \mu_j\}$ and $0 < \varepsilon ' < \varepsilon$.  
It follows that $\omega_i = \zeta_j$ for all $i \in \{1,\ldots, \mu_j\}$ and 
so  the polynomials $f(z)$ and $g(z)$ 
have the same roots with the same multiplicities. 
Therefore, $g(z) = cf(z)$ for some $c \neq 0$. 
This completes the proof. 
\end{proof}

The following result is an inverse to the theorem on continuity of roots.  

\bt
Let $\delta > 0$.  There exists $\varepsilon = \varepsilon(\delta) > 0$ such that 
if the polynomials 
$f(z)$ and $g(z)$ of degree $n$ are $\varepsilon$-aligned, then 
$g(z)$ is a $\delta$-deformation of $cf(z)$ for some $c \neq 0$. 
\et

\begin{proof}
This follows from Vi\' ete's formulae for the coefficients of a polynomial as 
the elementary symmetric functions of the roots. 
For 
$k \in \{1,2,\ldots, n\}$, the $k$th elementary symmetric function of $n$ variables is 
the polynomial 
\[
\sigma_k(z_1,\ldots, z_n) = \sum_{1 \leq i_1 < \cdots < i_k \leq n} z_{i_1} \cdots z_{i_k}. 
\]
We define $\sigma_0(z_1,\ldots, z_n) = 1$. 

Let 
\[
f(z) =  \sum_{i=0}^{n} a_i z^i = a_n \prod_{j=1}^n ( z-\zeta_j)    
\] 
with $a_n \neq 0$.    Vi\' ete's formulae are the identities 
\[
\frac{a_i}{a_n} = (-1)^{n-i} \sigma_{n-i}(\zeta_1,\ldots, \zeta_n)  
\]
for all $i \in \{0,1,\ldots, n\}$. 

Let $\varepsilon > 0$ and let $g(z)$ 
be a polynomial of degree $n$ that is $\varepsilon$-aligned to $f(z)$. 
We enumerate the roots $\omega_1,\ldots, \omega_n$ of $g(z)$ so that 
\[
 |\omega_j - \zeta_j| < \varepsilon 
 \]
 for all $j \in \{1,2,\ldots, n\}$. 
Setting 
\[
t_j   = \omega_j - \zeta_j  
\]
we have 
\[
 |t_j| < \varepsilon 
 \]
 for all $j \in \{1,2,\ldots, n\}$. 
There exists $b_n \neq 0$ such that  
\[
g(z)  = b_n \prod_{j=1}^n ( z-\omega_j)  = b_n \prod_{j=1}^n ( z - \zeta_j - t_j)  =  \sum_{i=0}^{n} b_i z^i 
\]
where  
\[
\frac{b_i}{b_n} = (-1)^{n-i} \sigma_{n-i}(\omega_1,\ldots, \omega_n)  
= (-1)^{n-i}  \sigma_{n-i}( \zeta_1 + t_1, \ldots,  \zeta_n + t_n)  
\]
for all $i \in \{0,1,\ldots, n-1\}$. 

The symmetric functions are polynomials, hence continuous, 
and so there exists $\varepsilon = \varepsilon(\delta) > 0$  
such that if $|t_j| < \varepsilon$ for all $j \in \{1,\ldots, n\}$, then 
\[
\left| \sigma_{n-i}( \zeta_1 + t_1, \ldots,  \zeta_n + t_n) - \sigma_{n-i}( \zeta_1, \ldots,  \zeta_n) \right| 
< \frac{\delta}{|b_n|}. 
\]
 Equivalently, 
\[
\left|  \frac{b_i}{b_n} - \frac{a_i}{a_n} \right|  < \frac{\delta}{|b_n|}
\]
and
\[
\left| b_i - c a_i  \right|  < \delta
\]
for $c = b_n/a_n \neq 0$  and all $i \in \{0,1,\ldots, n\}$. 
Thus,  $g(z)$ is a $\delta$-deformation of $cf(z)$. 
This completes the proof. 
\end{proof}

\section{Polynomials with only the zero root}                        \label{ContRoot:section:no-nonzero} 
In this section we prove the continuity of roots of polynomials whose only root is zero, that is, 
polynomials of the form $f(z) = a_nz^n$ with $a_n \neq 0$.  

\bt            \label{ContRoot:theorem:NoNonzeroRoot} 
Let $f(z) = a_nz^n$ with $a_n \neq 0$.  
For all $\varepsilon > 0$ there exists $\delta > 0$ 
such that if the polynomial 
\[
g(z) = \sum_{i=0}^n b_i z^i
\]
is a $\delta$-deformation of $f(z)$, then $|\omega| < \varepsilon $ for all roots $\omega$ of $g(z)$, 
that is, $g(z)$ is $\varepsilon$-aligned to $f(z)$. 
\et

\begin{proof}
Let $0 < \varepsilon < 1$ and  let  
\beq                           \label{ContRoot:delta}
0 < \delta <   \frac{\varepsilon^{n} |a_n|}{2n}. 
\eeq 
If $g(z) = \sum_{i=0}^n b_i z^i$ is a $\delta$-deformation of $f(z)$, then 
\[
|b_n - a_n| < \delta 
\]
and 
\[
|b_i| < \delta \qquad \text{for all $i \in \{0,1,\ldots, n-1\}$.} 
\]
The inequality 
\[
|a_n| - |b_n| \leq |b_n - a_n| < \delta  <  \frac{ |a_n|}{2n}
\]
implies 
\[
|b_n| >  |a_n| - \frac{ |a_n|}{2n} \geq \frac{ |a_n|}{2} > 0.
\]
Therefore,
\beq                            \label{ContRoot:bibn}
\frac{|b_i|}{|b_n|} < \frac{\delta}{|a_n|/2} =  \frac{2\delta}{|a_n|} 
\eeq
for all $i \in \{0,1,\ldots, n-1\}$.

If $|\omega| \geq 1$ and 
\[
g(\omega) = \sum_{i=0}^n b_i \omega^i = 0
\]
then 
\begin{align*}
0 < \frac{|a_n|}{2}  |\omega|^n  & < |b_n| |\omega|^n = \left| -b_n\omega^n \right| = \left|  \sum_{i=0}^{n-1} b_i \omega^i \right| \\
& 
\leq \sum_{i=0}^{n-1} \left|  b_i \right|   \left|  \omega\right|^i 
< \delta \sum_{i=0}^{n-1}   \left|  \omega\right|^i  \leq n\delta |\omega|^{n-1} 
\end{align*}
and so 
\[
|\omega| < \frac{2n\delta}{|a_n|} < \varepsilon^n < 1
\]
which is absurd.  Therefore, if $\omega \neq 0$ and $g(\omega) = 0$, 
then $0 < |\omega| < 1$. 

If $0 < |\omega| < 1$ and  $g(\omega) = b_n \omega^n + \sum_{i=0}^{n-1} b_i \omega^i = 0$, then 
\[
-1 = \sum_{i=0}^{n-1} \frac{b_i\omega^i}{b_n\omega^n} 
= \sum_{i=0}^{n-1} \frac{b_i}{b_n} \frac{1}{\omega^{n-i}}. 
\]
Applying inequalities~\eqref{ContRoot:delta} and~\eqref{ContRoot:bibn}, we obtain  
\begin{align*}
1 &\leq  \sum_{i=0}^{n-1} \frac{|b_i|}{|b_n|} \frac{1}{|\omega|^{n-i}} 
 <  \frac{2\delta}{ |a_n|} \sum_{i=0}^{n-1}\frac{1}{|\omega|^{n-i}} \\
&  \leq  \left(  \frac{2n\delta}{ |a_n|} \right)  \frac{1}{|\omega|^{n}} 
< \left( \frac{\varepsilon}{|\omega|} \right)^{n}
\end{align*}
and so $|\omega| < \varepsilon$ for all roots $\omega$ of $g(z)$. 
This is the definition of ``continuity of roots'' for polynomials $f(z)$ whose only root is 0. 
This completes the proof. 
\end{proof}

\section{Polynomials with nonzero roots}                      \label{ContRoot:section:nonzero} 

We begin with a simple result that implies continuity of roots of a polynomial 
with only  simple roots, 
that is, roots  of multiplicity one. 

\bt                                              \label{ContRoot:theorem:AllRoots} 
Let $f(z)$ be a polynomial of degree $n$.  
For every $\varepsilon > 0$, there exists $\delta  = \delta(f,\varepsilon)f) > 0$ such that 
if the polynomial $g(z)$ is a $\delta$-deformation of $f(z)$, 
then for every root $\zeta$ of $f(z)$ 
there is a root $\omega$ of $g(z)$ such that $|\omega - \zeta| < \varepsilon$. 
\et

\begin{proof}                \label{ContRoot:theorem:AllRoots} 
Let 
\[
f(z) =  \sum_{i=0}^{n} a_iz^i \in K[z]
\]
be a  polynomial of degree $n$ with  distinct roots $\zeta_1,\ldots, \zeta_{\ell}$ 
and let 
\[
M = \max(1,|\zeta_1|,\ldots, |\zeta_{\ell}|).
\]  
For  $0 < \varepsilon  < M$, let  
 \beq                                          \label{ContRoot:AllRoots-delta-0} 
0 < \delta <    \frac{|a_n|}{2(n+1)} \left( \frac{\varepsilon}{M} \right) ^n.
\eeq
Let 
\[
g(z) =   \sum_{i=0}^{n}b_iz^i \in K[z]  
\]
be a  polynomial of degree $n$ that is a $\delta$-deformation of $f(z)$.  
The inequality 
\[
|b_n - a_n| < \delta < \frac{|a_n|}{2}
\]
implies 
\[
 \frac{|a_n|}{2} < |b_n|. 
\]
Let 
\[
g(z) = b_n \prod_{j=1}^n (z- \omega_j) 
\]
where $\omega_1, \ldots, \omega_n$ are the $n$ not necessarily distinct roots of  $g(z)$.  
For all $k \in \{1,\ldots, \ell\}$ we have $f(\zeta_k)=0$ and so 
\[
b_n \prod_{j=1}^{n} (\zeta_k - \omega_j) =  g(\zeta_k) =  g(\zeta_k) - f(\zeta_k) = \sum_{i=0}^{n} (b_i-a_i)\zeta_k^i. 
\]
Therefore, 
\begin{align*}
 \frac{|a_n|}{2} \left( \min_{1\leq j \leq n}   |\zeta_k - \omega_j| \right)^n  
 & \leq  |b_n| \prod_{j=1}^{n} |\zeta_k - \omega_j| \\ 
&  \leq \sum_{i=0}^{n} |b_i-a_i| |\zeta_k|^i \\ 
& \leq (n+1)\delta M^n
\end{align*}
and so there exists $j_k \in \{1,\ldots, n\}$ such that 
\[
 |\zeta_k - \omega_{j_k}|^n  \leq \frac{2 (n+1)\delta M^n}{|a_n|}.
 \]
Equivalently,  
\[
 |\zeta_k - \omega_{j_k}| \leq 
 \left(\frac{2 (n+1)\delta }{|a_n|}\right)^{1/n} M < \varepsilon.
\]
This completes the proof.  
\end{proof}

\bt                                              \label{ContRoot:theorem:AllRootsSimple} 
Let $f(z)$ be a polynomial of degree $n$ with only simple roots.  
For every $\varepsilon > 0$, there exists $\delta  = \delta(f,\varepsilon) > 0$ such that 
if the polynomial $g(z)$ is a $\delta$-deformation of $f(z)$, 
then $g(z)$ has only simple roots 
and, for every root $\zeta$ of $f(z)$, 
there is a root $\omega$ of $g(z)$ such that $|\omega - \zeta| < \varepsilon$. 
\et

\begin{proof}
 This follows immediately from Theorem~\ref{ContRoot:theorem:AllRoots}. 
\end{proof}

\bl                                          \label{ContRoot:lemma:ab-hat} 
Let 
\[
f(z) =  \sum_{i=0}^{n} a_iz^i \in K[z]
\]
be a polynomial of degree $n$ and let $\zeta$ be a nonzero  root of $f(z)$.  
Then 
\[
\hat{f}(z) = \frac{f(z)}{z-\zeta} =  \sum_{i=0}^{n-1} \hat{a}_iz^i   
\]
is a polynomial in $K[z]$ of degree $n-1$ and  
\[
\hat{a}_i = -\frac{1}{\zeta^{i+1}}  \sum_{k=0}^i  a_k \zeta^k 
\]
for all $i \in \{0,1,\ldots, n-1\}$. 
\el

\begin{proof}
This is a straightforward calculation.  We have 
\begin{align*}
 \sum_{i=0}^{n} a_iz^i  & = f(z)  = (z - \zeta) \hat{f}(z)
  = (z - \zeta) \sum_{i=0}^{n-1} \hat{a}_iz^i  \\
& = - \zeta \hat{a}_0 + \sum_{i=1}^{n-1} \left( \hat{a}_{i-1}  -  \zeta\hat{a}_i \right) z^i  + \hat{a}_{n-1} z^{n} 
\end{align*}
and so 
\[
a_i = \begin{cases}
- \zeta \hat{a}_0 & \text{if $i = 0$} \\
 \hat{a}_{i-1} -  \zeta \hat{a}_i &  \text{if $i \in \{1,\ldots, n-1\}$} \\
 \hat{a}_{n-1}  & \text{if $i = n$.} 
\end{cases}
\]
Equivalently, 
\[
 \hat{a}_0   = -\frac{a_0}{\zeta} 
\]
and 
\[
 \hat{a}_{i}  =  -\frac{1}{\zeta} \left( a_i - \hat{a}_{i-1} \right) \qquad \text{for $i \in \{1,\ldots, n-1\}$.} 
\]
We obtain 
\[
\hat{a}_1 = -\frac{1}{\zeta} (a_1 - \hat{a}_0 )   
 = -\frac{1}{\zeta^2} \left( a_0 + a_1 \zeta \right).
\]
If $i \in \{1,2,\ldots, n-1\}$ and 
\[
\hat{a}_{i-1} = -\frac{1}{\zeta^{i}}  \sum_{k=0}^{i-1}  a_k \zeta^k 
\]
then 
\begin{align*}
 \hat{a}_i  & =  -\frac{1}{\zeta} \left( a_i - \hat{a}_{i-1} \right) 
 =  -\frac{1}{\zeta} \left( a_i + \frac{1}{\zeta^{i}}  \sum_{k=0}^{i-1}  a_k \zeta^k  \right) \\
& =  -\frac{1}{\zeta^{i+1}}  \sum_{k=0}^{i}  a_k \zeta^k .
\end{align*}
This completes the proof. 
\end{proof}

\bl                                     \label{ContRoot:lemma:ab-delta} 
Let 
\[
f(z) =   \sum_{i=0}^{n} a_iz^i \in K[z]
\]
and  
\[
g(z) =  \sum_{i=0}^{n } b_iz^i \in K[z]
\]
be  polynomials of degree $n$. 
Let $\zeta_1$ be a nonzero root of $f(z)$,  let $\omega_1$ be a nonzero root of $g(z)$, 
and let 
\[
\hat{f}(z) = \frac{f(z)}{z-\zeta_1} =  \sum_{i=0}^{n-1} \hat{a}_iz^i
\]
and 
\[
\hat{g}(z) = \frac{g(z)}{z-\omega_1} =  \sum_{i=0}^{n-1} \hat{b}_iz^i. 
\]
For every $\lambda > 0$  there exists $\kappa > 0$ such that if 
$g(z)$ is a $\kappa$-deformation of $f(z)$ and if
\[
|\omega_1 - \zeta_1| < \kappa
\]
then $\hat{g}(z)$ is a $\lambda$-deformation of $\hat{f}(z)$. 
\el

\begin{proof} 
Let $K^{\times} =  K\setminus \{ 0\}$. 
For all $i \in \{0,1,\ldots, n-1 \}$, the function  
\[
\varphi_i:K^{i+1} \times K^{\times} \rightarrow K
\]
defined by 
\[
\varphi_i(u_0,u_1,\ldots, u_i,v) = -\frac{1}{v^{i+1}} \sum_{k=0}^i  u_k v^k 
\] 
is continuous. Let  
\[
(a_0,a_1,\ldots, a_n, \zeta_1) \in K^{n+1} \times K^{\times}.
\]
There exists $\kappa > 0$ such that if 
\[
(u_0,u_1,\ldots, u_n,v) \in K^{n+1} \times K^{\times} 
\]
and if 
\[
|u_i - a_i| < \kappa 
\]
for all $i \in \{0,1,\ldots, n \}$ and 
\[
|v - \zeta_1| < \kappa
\]
then 
\[
\left| \varphi_i(u_0,u_1,\ldots, u_i,v) - \varphi_i(a_0,a_1,\ldots, a_i, \zeta_1)  \right| < \lambda
\]
for all $i \in \{0,1,\ldots, n-1 \}$.

By Lemma~\ref{ContRoot:lemma:ab-hat}, 
\[
 \varphi_i(a_0,a_1,\ldots, a_i, \zeta_1) = \hat{a}_i
\]
and 
\[
 \varphi_i(b_0,b_1,\ldots, b_i,\omega_1) = \hat{b}_i 
\] 
for all $i \in \{0,1,\ldots, n-1 \}$.  
If $g(z)$ is a $\kappa$-deformation of $f(z)$, then 
\[
|b_i - a_i| < \kappa 
\]
for all $i \in \{0,1,\ldots, n\}$.  The inequality 
\[
|\omega_1 - \zeta_1| < \kappa
\]
implies  
\[
|\hat{b}_i  - \hat{a}_i | = 
\left|  \varphi_i(b_0,b_1,\ldots, b_i,\omega_1) -  \varphi_i(a_0,a_1,\ldots, a_i, \zeta_1) \right| 
 < \lambda
\]
for all $i \in \{0,1,\ldots, n-1 \}$ 
and so $\hat{g}(z)$ is a $\lambda$-deformation of $\hat{f}(z)$. 
This completes the proof.  
\end{proof}

Here is the continuity of roots theorem for polynomials 
with a nonzero root.

\bt                              \label{ContRoot:theorem:ContinuityRoots} 
Let $f(z)$ be a  polynomial of degree $n \geq 1$ with a nonzero root.    
For all $\varepsilon > 0$, there exists $\delta =  \delta(f,\varepsilon) > 0$ 
such that if the polynomial $g(z)$ is a $\delta$-deformation of $f(z)$, then 
$g(z)$ is $\varepsilon$-aligned to $f(z)$. 
\et

\begin{proof}
Let 
\[
f(z) =  \sum_{i=0}^{n} a_iz^i \in K[z].
\] 
Let $\zeta_1,\ldots, \zeta_{\ell}$ be the distinct roots of $f(z)$ with $\zeta_1 \neq 0$.  
Without loss of generality we can assume that $0 < \varepsilon \leq |\zeta_1|$ 
and that $\varepsilon$ satisfies inequality~\eqref{ContRoot:epsilon}.  

The proof is by induction on the degree $n$ of $f(z)$. The case $n=1$ follows 
from Theorem~\ref{ContRoot:theorem:AllRoots}.  
Let $n \geq 2$ and assume that Theorem~\ref{ContRoot:theorem:ContinuityRoots} 
is true for  polynomials of degree $n-1$.  

 The polynomial 
\[
\hat{f}(z) = \frac{f(z)}{z-\zeta_1} = \sum_{i=0}^{n-1} \hat{a}_i z^i 
\]
has degree $n-1$.   
By the induction hypothesis, there exists $\lambda  = \lambda(\hat{f},\varepsilon) > 0$ such that 
if the polynomial $h(z)$ is a $\lambda$-deformation of  $\hat{f}(z)$, 
then $h(z)$ is $\varepsilon$-aligned to $\hat{f}(z)$. 
It suffices to prove that there exists $\delta   =  \delta(f,\varepsilon) > 0$ such that if the  polynomial 
\[
g(z) =  \sum_{i=0}^{n} b_iz^i \in K[z] 
\] 
is a $\delta$-deformation of $f(z)$, then $g(z)$ has a nonzero root $\omega_1$ with 
$|\omega_1 - \zeta_1| < \varepsilon$ such that the  polynomial of degree $n-1$ 
\[
\hat{g}(z) = \frac{g(z)}{z- \omega_1} = \sum_{i=0}^{n-1} \hat{b}_i z^i 
\]  
is a $\lambda$-deformation of $\hat{f}(z) $. 

 Let $0 < \kappa < \varepsilon$. 
By Theorem~\ref{ContRoot:theorem:AllRoots}, there exists $\delta_1 = \delta_1(f,\kappa) > 0$ such that 
if 
\[
g(z) =   \sum_{i=0}^{n} b_iz^i  
\]
is a   polynomial of degree $n$ that is a $\delta_1$-deformation of $f(z)$, 
then there is a root $\omega_1$ of $g(z)$ such that 
\[
|\omega_1 - \zeta_1| < \kappa < \varepsilon.  
\]
The inequality  $\varepsilon \leq |\zeta_1|$ implies $\omega_1 \neq 0$. 

By Lemma~\ref{ContRoot:lemma:ab-delta}, there exists $0 < \kappa  < \varepsilon$ 
such that if the polynomial $g(z)$ 
is a $\kappa$-deformation of $f(z)$ and if $\omega_1$ is a nonzero root of $g(z)$ 
with $|\omega_1 - \zeta_1| < \kappa$, then 
$\hat{g}(z)$ is a  $\lambda$-deformation of $\hat{f}(z)$. 

Let $\delta = \min(\kappa, \delta_1)$.  If $g(z)$ is a $\delta$-deformation of $f(z)$, 
then $g(z)$ has a nonzero root $\omega_1$ with $|\omega_1 - \zeta_1| < \kappa$ 
such that $\hat{g}(z)$ is a  $\lambda$-deformation of $\hat{f}(z)$.  
It follows that $\hat{g}(z)$ is a $\varepsilon$-aligned to  $\hat{f}(z)$.

Let $\zeta_1,  \zeta_2, \ldots, \zeta_n$ be the roots of $f(z)$ 
and let $\omega_1,\omega_2,\ldots, \omega_n$ 
be the roots of $g(z)$.  
We have 
\[
\hat{f}(z) = \frac{f(z)}{z- \zeta_1} = a_n \prod_{j=2}^n (z-\zeta_j) 
\]
and
\[
\hat{g}(z) = \frac{g(z)}{z- \omega_1} = b_n \prod_{j=2}^n (z-\omega_j). 
\]
The polynomials $\hat{f}(z)$ and $\hat{g}(z)$ are $\varepsilon$-aligned.  
Because 
\[
|\omega_1 - \zeta_1| < \kappa < \varepsilon  
\]
it follows that the polynomals $f(z)$ and $g(z)$ are $\varepsilon$-aligned. 
This completes the proof. 
\end{proof}

\end{document}